\documentstyle[11pt]{article}

\newtheorem{th}{Theorem}[section]
\newtheorem{de}[th]{Definition}

\newtheorem{re}[th]{Remark}

{

\begin{document}

\large
\title {\bf GLOBAL VERSION OF BJORNESTAL'S ESTIMATE FOR METRIC
PROJECTION OPERATOR IN BANACH SPACE}
\author{\bf Ya. I. Alber $^{*}$\\
Department of Mathematics\\
Technioin-Israel Institute of Technology\\
Haifa 32000, Israel}
\date{}
\maketitle

{\bf Abstract} $\;$- In 1979, B.Bjornestal obtained local estimate for a
modulus of uniform continuity of metric projection operator on closed
subspace in uniformly convex and uniformly  smooth Banach space $B$. In the
present paper we give the global version of this result for the  projection
operator on an arbitrary closed convex set in $B$.

\section{Introduction and Preliminaries}
\setcounter{equation}{0}

Metric projection operators $P_\Omega$ on convex closed sets $\Omega$
(in the sense of best approximation)
are widely used in theoretical and applied  areas of
mathematics especially connected with problems  of optimization  and
approximation.
As examples one can provide iterative-projection methods for solving
equations, variational inequalities and minimization of functionals \cite{a1}
and  methods of alternating projections for finding common point of convex
closed sets in Hilbert spaces \cite{d,bc,bf}.

Let us remind a definition of the metric projection operator.
Let $B$ be a real, uniformly convex and uniformly  smooth  (reflexive)
Banach space with $B^*$ its dual space,
$\Omega$ be  a closed convex set in $B,\; <w,v>$ a dual product
in $B$, i.e. pairing  between
$w \in  B^*$ and $v \in  B \; ((y,x)$ is inner product in Hilbert
space $H$, if we identify $H$ and $H^*$);
the signs $||\cdot||$ and $ \;||\cdot||_{B^{*}}$  denote
the norms in Banach spaces $B$ and $B^{*}$, respectively.
\begin {de}  \label{gpo}
The operator $P_{\Omega}: B \rightarrow {\Omega} \subset B$
is called metric projection operator if it yields
the correspondence between an arbitrary point  $x \in  B$ and nearest
point $\bar x \in \Omega $ according to minimization problem  \\
 \end {de}
------------------------\\
$^{*}$ This research was supported in part by the Ministry of Science Grant
3481-1-91 and by the Ministry of Absorption Center for Absorption in Science.

\begin{equation} \label{f1}
P_{\Omega}x = \bar x; \;\;\; \bar  x: ||x - \bar x|| =
\inf_{\xi \in \Omega} ||x - \xi||.
\end{equation}

Metric projection operators have extremely good properties in Hilbert spaces
\cite{h,z,a1}. In the contrary,
they are too poor in Banach spaces. For instance,
operators $P_\Omega$ do not possess such important properties as monotonicity,
non-expansiveness, absolutely best approximation \cite{a1},
which make the metric projection operators in Hilbert spaces
exceptionally effective ones. To illustrate this let us remind the properties
of metric projection operator on subspace $M$ of Hilbert space. Here, such
operator is is  orthogonal, linear, non-expansive, self-adjoint and idempotent
\cite{d}. Metric projection operator on subspace $M$  of Banach space has no
these properties  in general \cite{fl}.

However,  $P_\Omega$ possess a number of good qualities realized in very
important applications \cite{a2,a3}. So, it is  uniformly continuous in
Banach space on each bounded set
and satisfies  basic variational principle \cite{li}
\begin{equation} \label{b1}
< J(x - \bar x), \bar x -\xi > \ge 0,  {\qquad }
\forall \xi \in\Omega,
\end{equation}
Here $J: B \rightarrow B^* $ a normalized duality mapping in $B$ \cite{li,a1}.

The smoothness properties  of metric projection operator are being studied
for a long time.  In Hilbert space it satisfies the Lipschitz
condition and, consequently, it is uniformly continuous. It is known that in
uniformly convex Banach space metric projection operator is always
continuous but not always uniformly continuous.

The results of F.Murray and J.Lindenstrauss (see \cite{h}) put forward
the problem:
"Whether or not the operator $P_\Omega$ in  uniformly convex and uniformly
smooth Banach space  is  uniformly continuous?" In 1979, B.Bjornestal
obtained positive answer to this question in form of estimate \cite{bj}
\begin {equation} \label{b2}
||P_M x - P_M y|| \leq 2\delta^{-1}_B (2\rho _B {(6||x - y||)}),
\end{equation}
where $M$ is a closed linear subspace of $B$,
$\rho _B (\tau)$ is  a modulus of smoothness,
$\delta _B (\epsilon) $  is a modulus of convexity of the space B,
and  $\delta^{-1}_B (\cdot)$
is the inverse function to  $\delta _B (\epsilon) $ \cite{lz}. {\it But
this result was only local \/} (it is fulfilled  if $x$ and $y$ are
sufficiently near  with each other  and $||x - \bar x|| = 1,$
$||y - \bar y|| = 1$).

Recently,  in paper \cite{zr} the following global estimate was
established in uniformly convex and uniformly smooth Banach space $B$
\begin {equation} \label{b3}
||P_\Omega x - P_\Omega y|| \leq ||x - y|| + 4 C_1 \delta^{-1}_B (N \psi
(||x - y|| / C_1))
\end{equation}
where $N$ is some fixed constant, $C_1 =  ||x - P_\Omega y||\bigvee
||P_\Omega x - y||,$ and $\psi$ is the function defined by the formula
$$\psi (t) = \int_{0}^{t} {\frac {\rho _B (s)}{s}} ds .$$

In \cite{an1} we obtained  another estimate of uniform continuty
of metric projection operator
in  uniformly convex and uniformly smooth Banach space $B$
\begin{equation} \label{b4}
||P_\Omega x - P_\Omega y|| \leq  C g^{-1}_B (NC g^{-1}_{B^{*}} {(N||x -
y||)}),
\end{equation}
where
$g_B (\epsilon)  = \delta  _B (\epsilon)  / \epsilon$,  $g^{-1}_B (\cdot)$
is an inverse function,
$N = 2LC$, $L$ is  constant, $1 < L  < 3.18,$ (see \cite{fg}) and
$$ C = 2  {max} \lbrace 1, ||x - P_\Omega y ||, ||y - P_\Omega x||\rbrace .$$
The estimate (\ref {b4}) looks more preferable because of the
integral in  (\ref {b3}).

However, the simple calculations show, that Bjornestal's estimate (\ref {b2})
better
than (\ref {b4}), first of all, comparing their orders. For instance,
known estimates  for the moduli of convexity and smoothness
of the spaces $l^p, L^p $  and $ W^p_m$, where  $\infty > p >1, $
$$ \rho_B (\tau) \leq p^{-1} \tau^p,\;\;
\delta _{B}(\epsilon) \ge (p-1)/8 \epsilon^2,\;\;     1 < p \leq 2 ,$$
$$\rho_B (\tau) \leq (p-1) \tau^2,\;\;
\delta _{B}(\epsilon) \ge p^{-1}{(\epsilon/2)}^p,\;\;    \infty > p > 2 $$
give the following orders: for (\ref {b2})
$$||P_\Omega x - P_\Omega y|| \sim ||x-y||^{2/p}, \;\; p \ge 2  ,$$
and for (\ref {b4})
$$||P_\Omega x - P_\Omega y|| \sim ||x-y||^{1/(p-1)}, \;\; p \ge 2  $$
Our remark for $p > 2$ is valid because $2/p > 1/(p-1)$. \\
Let now $1 < p \leq 2.$ Then the estimate (\ref {b2}) yields
$$||P_\Omega x - P_\Omega y|| \sim ||x-y||^{p/2}, \;\;  2 \ge p > 1,  $$
and for  estimate (\ref {b4}) we have
$$||P_\Omega x - P_\Omega y|| \sim ||x-y||^{p-1}, \;\; 2 \ge p > 1   .$$
Our remark for $1 < p < 2$ is valid because $p/2 > p-1 $. For $p = 2$
(Hilbert case) (\ref {b2}) and (\ref {b4}) give the same orders:
$$||P_\Omega x - P_\Omega y|| \sim ||x-y|| .$$

Let us emphasize again that (\ref {b2}) is local estimate. In the next
Section we will obtain its global form for the arbitrary closed convex
set $\Omega$ in Banach spaces.

\section{Auxiliary theorems}
\setcounter{equation}{0}

The lower and upper parallelogram inequalities and
estimates of duality mappings in uniformly convex and uniformly smooth Banach
spaces (respectively) obtained first in \cite{an2,an3,n}
are used as the basis in order to prove uniform  continuity of the metric
projection  operators in Banach spaces.
In this Section we will have proven two auxiliary theorems.
\begin {th} \label{l1}
In uniformly smooth Banach space $B$ the following estimate
\begin{equation} \label{d1}
<Jx -Jy, x - y> \leq 8||x - y||^2 + C \rho _B {(||x - y||)}, {\qquad }
\forall x,y \in B
\end{equation}
is valid, where
$$ C = C(||x||,||y||) = 4\; {max} \lbrace 2L, ||x|| + ||y|| \rbrace .$$
\end {th}
{\bf Proof}.  Denote
$$D = 2^{-1}( ||x||^2 + ||y||^2 - ||2^{-1}(x + y)||^2 /2)$$
and consider two possibilities: \\
\\
(i) Let  $||x + y|| \leq ||x - y||$. Then
$$||x|| + ||y|| \leq ||x + y|| + ||x - y|| \leq 2||x - y|| .$$
Involving this expression in square, we obtain
$$ 2^{-1}||x||^2 + 2^{-1}||y||^2 + ||x||||y|| \leq 2 ||x - y||^2 .$$
Now, let us subtract $||2^{-1}(x + y)||^2$ from both parts of this
inequality. We have
$$ D \leq  2 ||x - y||^2  - (||2^{-1}(x + y)||^2 + ||x|||| y||). $$
If $||2^{-1}(x + y)||^2 + ||x|||| y|| \ge  ||x - y||^2 $ then immediately
\begin{equation} \label{b5}
D \leq ||x - y||^2.
\end{equation}
Suppose that the contrary inequality occurs. In this case, it is easily
verify
$$ 2^{-1}||x||^2 + 2^{-1}||y||^2 - ||2^{-1}(x + y)||^2  $$
$$\leq ||2^{-1}(x + y)||^2 + ||x||||y|| \leq ||x - y||^2$$
following from the estimate
$(||x|| - ||y||)^2  \leq ||x + y||^2.$
i.e. (\ref{b5}) is valid. \\
\\
(ii) Let now  $||x + y|| \ge ||x - y||$. It can be shown that
\begin{equation} \label{b6}
||x|| + ||y|| - ||x + y||  \leq \epsilon (x,y)
\end{equation}
where
\begin{equation} \label{b7}
\epsilon (x,y) = ||x + y|| \; \rho _B ({\frac {||x - y||}{||x + y||}}).
\end{equation}
Indeed, let us replace
$$x = {\frac 12} (u + v),\;\;y =  {\frac 12} (u - v)$$
and set
$$\alpha = {\frac {u}{||u||}},\;\;\beta = {\frac {v}{||u||}} .$$
Using the definition of the modulus of smoothness $ \rho_B (\tau)$,
one can write
$$||x|| + ||y|| - ||x + y||  = 2^{-1}(||u+v|| + ||u-v||) - ||u|| $$
$$ = 2^{-1}||u||(||\alpha + \beta|| +  ||\alpha - \beta|| -2)  $$
$$ \leq ||u||\;{sup}[2^{-1}(||\alpha + \beta|| +  ||\alpha - \beta||) -1,\;
||\alpha || = 1,\; ||\beta|| = \tau ] $$
$$ \leq ||u|| \rho_B (||\beta ||) .$$
Returning to old notations we obtain (\ref {b6}) and (\ref {b7}). Thus,
$$ ||{\frac {x + y}{2}}|| \ge {\frac {||x|| + ||y|| - \epsilon (x,y)}{2}} .$$
It's right hand part is nonnegative. In fact, using the property
$ \rho_B (\tau) \leq \tau $ \cite{lz} we establish inequality
$$ ||x|| + ||y|| - \epsilon (x,y) \ge ||x|| + ||y|| - ||x-y|| \ge 0 .$$
Then
$$  ||{\frac {x + y}{2}}||^2 \ge
({\frac {||x|| + ||y||}{2}})^2 -  \epsilon (x,y) {\frac {||x|| + ||y||}{2}}$$
In virtue of $||x|| - ||y|| \leq ||x-y||$ one have
\begin{equation} \label{b8}
D \leq ({\frac {||x|| - ||y||}{2}})^2
+  \epsilon (x,y) {\frac {||x|| + ||y||}{2}} \leq
||{\frac {x - y}{2}}||^2 +  \epsilon (x,y) {\frac {||x|| + ||y||}{2}}.
\end{equation}

a) Suppose that $ ||x + y|| < 1 $ then $ ||x + y||^{-1} ||x - y|| > ||x - y||.$
It is known \cite{fg} that the inequality
\begin{equation} \label{b9}
\tau{_1}^2 \rho_B (\tau_2) \leq L \tau{_2}^2 \rho_B (\tau_1),\;\;
0 \leq \tau_1 \leq \tau_2,\;\;1 < L < 3.18
\end{equation}
takes place in arbitrary Banach space. By (\ref {b8}) and (\ref {b9})
$$ \rho_B ( ||x - y|| / ||x + y||) \leq L ||x + y||^{-2} \rho_B (||x - y||).$$
It follows from the last estimate that
$$D \leq 4^{-1} ||x - y||^2 + 2^{-1}L (||x|| + ||y||)
||x + y||^{-1} \rho_B (||x - y||) .$$
So far as  $||x + y|| \ge ||x - y||$  we have
$$ 2^{-1}||x + y||^{-1}(||x|| + ||y||) \leq
(2 ||x + y||)^{-1} (||x + y|| + ||x - y||) \leq 1 .$$
Therefore
\begin{equation} \label{b10}
D \leq 4^{-1} ||x - y||^2 + L \rho_B (||x - y||).
\end{equation}b) Let us now  $ ||x + y|| \ge 1, $ Then we obtain an addition
to
(\ref {b10}) in the form
\begin{equation} \label{b11}
D \leq 4^{-1} ||x - y||^2 + 2^{-1}(||x|| + ||y||) \rho_B (||x - y||).
\end{equation}
Here we accounted (\ref {b8})   and the convexity of $\rho_B (\tau)$.
The estimates (\ref {b5}), (\ref {b10}) and (\ref {b11}) being joined
together give
$$2||x||^2 + 2||y||^2 + ||x + y||^2 \leq 4 ||x - y||^2 + 2 {max}
\lbrace 2L, ||x|| + ||y|| \rbrace \rho_B (||x - y||) .$$
This is {\it the upper parallelogram inequality \/} in uniformly smooth
Banach space \cite{an3}.

Denote the right hand part of this inequality by $k(||x -y||)$. Then
\begin{equation} \label{b11}
D \leq  k(||x -y||)/4.
\end{equation}

Further, for convex function  $\phi(x) = ||x||^2 / 2 $, let us construct
the concave (with respect to $\lambda$) function
$$\Phi (\lambda) = \lambda \phi(x) + (1 - \lambda) \phi(y) - \phi(y +
\lambda (x -y)),\;\;0 \leq \lambda \leq 1 .$$
It is obvious that $\Phi (0) = 0.$  Suppose  $0 < \lambda_1 \leq
\lambda_2$. Then
$ {\lambda_1}^{-1} \Phi (\lambda_1) \ge {\lambda_2}^{-1} \Phi (\lambda_2),$
i.e. $(\Phi (\lambda) / \lambda)' \leq 0.$
Opening this expression we have $\Phi'(\lambda) \leq \Phi (\lambda)/ \lambda.$
In particular,  $\Phi'(1/4) \leq  4\Phi (1/4).$  But
$$\Phi (1/4) = {\frac {1}{4}} \phi(x) +  {\frac {3}{4}} \phi(y) -
\phi({\frac {1}{4}}x + {\frac {3}{4}}y) .$$
It follows from (\ref {b11}) that for all $z_1$ and $z_2$
$$\phi({\frac {z_1 + z_2}{2}}) \ge  {\frac {\phi (z_1)}{2}} +
{\frac {\phi (z_2)}{2}} - k(||x -y||)/8 .$$
Let us set $z_1 = (x+y)/2$ and $z_2 = y$ and use the property
$k(t/2) \leq k(t)/2.$ We obtain
$$\phi({\frac {1}{4}}x + {\frac {3}{4}}y) =
\phi({\frac {1}{2}}({\frac {1}{2}}x + {\frac {1}{2}}y) + {\frac {1}{2}}y)
\ge {\frac {1}{2}}\phi({\frac {x + y}{2}}) + {\frac {1}{2}} \phi(y) -
{\frac {1}{8}} k(||{\frac {x -y}{2}}||) $$
$$\ge {\frac {1}{4}} \phi(x) + {\frac {1}{4}} \phi(y) -
{\frac {1}{16}} k(||x -y||) + {\frac {1}{2}} \phi(y) -  {\frac {1}{8}}
k(||{\frac {x - y}{2}}||)  $$
$$ = {\frac {1}{4}} \phi(x) + {\frac {3}{4}} \phi(y) -
{\frac {1}{8}}k(||x -y||) .$$
Thus, $ \Phi (1/4) \leq k(||x -y||)/8 $ and
$$\Phi'(1/4) =  \phi(x) - \phi(y) - <\phi'(y +
{\frac {1}{4}}(x -y)), x-y> \leq k(||x -y||)/2 .$$
Further, write down this inequality with $y$ and $x$ in places of $x$ and $y:$
$$ \phi(y) - \phi(x) - <\phi'(x +
{\frac {1}{4}}(y -x)), y-x> \leq k(||x -y||)/2 .$$
The two last inequalities  give  in sum
$$<\phi'(x + {\frac {1}{4}}(x -y)) - \phi'(y + {\frac {1}{4}}(x -y)),
x-y> \leq k(||x -y||) .$$
One can make now non-degenerate substitution of the variables $x$ and $y$
$$ z_1 = 2x - {\frac {1}{2}}(x -y),\;\; z_2 = 2y + {\frac {1}{2}}(x -y) $$
which leads to relations
$$z_1 - z_2 = x - y\;\; {and}\;\; ||x|| + ||y|| \leq ||z_1|| + ||z_2|| .$$
Taking into consideration the fact, that $Jx = \phi'(x)$ is homogeneous
operator, we find
$$  <Jz_1 -  Jz_2,z_1 - z_2> \leq 2 k(||z_1 - z_2||) .$$
Theorem is completely proved.

The following inequality is proved shorter than the previous
but it has the constant $L$ and
the function $C(||x||,||y||)$ under the sign of the modulus of smoothness
$\rho _B (\tau).$
\begin {th} \label{l16}
In uniformly smooth Banach space $B$ the estimate
\begin{equation} \label{d1}
<Jx -Jy, x - y> \leq (2L)^{-1} \rho _B {(8CL||x - y||)}, {\qquad }
\forall x,y \in B
\end{equation}
is valid, where
$$ C = C(||x||,||y||) = 2 {max} \lbrace 1, \sqrt{(||x||^2 + ||y||^2)/2}
\rbrace $$
\end {th}
{\bf Proof}.  Lemma 2.1 from \cite{an1} (cf. Theorem 2 from \cite{an2}) gives
the following estimate
\begin{equation} \label{a7}
<Jx - Jy, x - y> \ge (2L)^{-1} \delta_{B^{*}} (||Jx-Jy||_{B^{*}}/C)
\end{equation}
for uniformly smooth space $B.$  From (\ref{a7}) we have
$$||Jx-Jy||_{B^{*}}||x - y|| \ge (2L)^{-1} \delta_{B^{*}}
(||Jx-Jy||_{B^{*}}/C) .$$
Since $g_{B^{*}} (\epsilon)  = \delta_{B^{*}} (\epsilon)  / \epsilon $,
we can write
\begin{equation} \label{a8}
g_{B^{*}} (||Jx-Jy||_{B^{*}}/C) \leq 2CL ||x - y||.
\end{equation}
It is known from geometry of Banach spaces  \cite{lz}  that
$$\rho _B (\tau) \ge \epsilon \rho /2 - \delta _{B^{*}} (\epsilon),\;\;
0 \leq \epsilon \leq 2,\;\; \tau > 0 .$$
Therefore
$$\rho _B (4\delta _{B^{*}}(\epsilon)/ \epsilon) \ge
\delta _{B^{*}}(\epsilon) .$$
We denote $h _B (\tau) = \rho _B (\tau) /\tau $. Then
$$h _B (4g_{B^{*}} (\epsilon)) \ge \epsilon / 4 .$$
Setting
$$\epsilon = ||Jx-Jy||_{B^{*}}/C, $$
and using non-decreasing of the function $h _B (\tau)$, we find from
(\ref{a8})
$$h _B (4g_{B^{*}} (\epsilon)) \leq h_B (8CL||x - y||) .$$
Therefore
\begin{equation} \label{p1}
||Jx-Jy||_{B^{*}} \leq 4C h_B (8CL||x - y||).
\end{equation}
Now, (\ref {d1}) can be obtained from the inequality of Cauchy-Schwarz.
Theorem is proved.
\begin{re}
$C(||x||,||y||)$ in estimates
(\ref{l1}) and (\ref{d1}) are absolute constants $C = 8 {max} \lbrace L,
R \rbrace $  and $C = 2 {max} \lbrace 1, R \rbrace $,  respectively
if $||x|| \leq R$ and $||y|| \leq R$. In these
cases,  (\ref{l1}) and (\ref{d1}) are quantitative description of the
property of uniform continuity (in the form of dual product)
for normalized duality mapping $J$.  At the same time (\ref{p1}) gives the
modulus of uniform continuity of $J$ in traditional form.
\end{re}

\section{Main theorems}
\setcounter{equation}{0}

In this Section we will have provided the estimates of continuity for metric
projection operator on convex closed set of
uniformly convex and uniformly  smooth Banach space $B$.  In the case
when $||x||$ and $||y||$ are uniformly bounded,
they are the estimates of moduli of uniform continuity of this operator on
each bounded set of $B$.
\begin {th} \label{th1}
In uniformly convex and uniformly smooth Banach space $B$ the following
estimate
\begin{equation} \label{b14}
||P_\Omega x - P_\Omega y|| \leq C \delta^{-1}_B (2LC_1
\rho _B {(||x - y||)}),
\end{equation}
is satisfies where$$C =  2 {max} \lbrace 1, ||x - \bar y||, ||y - \bar x||
\rbrace,$$
$$ C_1 =  16 + 24\; {max} \lbrace L, ||x - \bar y||, ||y - \bar x||
\rbrace $$
\end{th}
{\bf Proof}. It is known \cite{an2} that
$$ \rho _B {(\tau)} \ge  \rho _H {(\tau)} =  \sqrt{1 + \tau^2} - 1 \ge
\tau^2 / (\tau + 2) .$$
Then
$$  \rho _B {(||x - y||)} \ge  ||x - y||^2 / (||x - y|| + 2) .$$
{}From this we have
$$||x - y||^2 \leq (||x|| + ||y|| + 2) \rho _B {(||x - y||)} ,$$
and taking (\ref {l1}) in consideration, we obtain
$$<Jx -Jy, x - y> \leq \bar {C_1} \rho _B {(||x - y||)} $$
where
$$ \bar {C_1} = 16 + 24\; {max} \lbrace L, ||x||, ||y|| \rbrace .$$
Now let us turn to making  the estimates of convex functional
$\varphi (x) = ||x||^2.$ We have
$$||x - \bar y||^2 - ||y - \bar y||^2 \leq 2 <J(y - \bar y), x - y>  +
<J(x - \bar y) - J(y - \bar y) , x - y>  $$
$$\leq 2 <J(y - \bar y), x - y>  +  C_2 \rho _B {(||x - y||)}, $$
$$ C_2 =  16 + 24\;{max} \lbrace L, ||x - \bar y||, ||y - \bar y|| \rbrace .$$
By analogy with previous, we can write
$$||y - \bar x||^2 - ||x - \bar x||^2 \leq 2 <J(x - \bar x), y - x>  +
C_3 \rho _B {(||x - y||)}, $$
$$ C_3 = 16 + 24\;{max} \lbrace 2L, ||x - \bar x||, ||y - \bar x|| \rbrace .$$
Add the last two inequalities. Then
$$(||x - \bar y||^2 - ||x - \bar x||^2) + (||y - \bar x||^2 -
||y - \bar y||^2)  $$
$$\leq 2 <J(y - \bar y) - J(x - \bar x) , x - y> + 2C_4 \rho _B {(||x - y||)},
$
$$ C_4 =  16 + 24\; {max} \lbrace L, ||x - \bar y||, ||y - \bar x||
\rbrace .$$
It is known \cite{li} that
$$<J(x - \bar x) - J(y - \bar y) , x - y> \ge 0 .$$
Therefore
\begin{equation} \label{b12}
(||x - \bar y||^2 - ||x - \bar x||^2) + (||y - \bar x||^2 -
||y - \bar y||^2 \leq  2C_4 \rho _B {(||x - y||)}.
\end{equation}
The conditions of uniform convexity of the functional
$\varphi (x) = ||x||^2 $ and uniform convexity of the Banach space $B$ gives
$$(||x - \bar y||^2 - ||x - \bar x||^2) \ge 2<J(x - \bar x),
 \bar x - \bar y> + (2L)^{-1} \delta_B (||\bar x - \bar y||/C_5),$$
$$ C_5 = 2 {max} \lbrace 1, ||x - \bar y||, ||x - \bar x|| \rbrace $$
and
$$(||y - \bar x||^2 - ||y - \bar y||^2) \ge 2<J(y - \bar y),
 \bar x - \bar y> + (2L)^{-1} \delta_B (||\bar x - \bar y||/C_6),$$
$$ C_6 = 2 {max} \lbrace 1, ||y - \bar x||, ||y - \bar y|| \rbrace .$$
Using
$$<J(x - \bar x) - J(y - \bar y) , \bar x - \bar y> \ge 0 ,$$
\cite{li}, we obtain
\begin{equation} \label{b13}
(||x - \bar y||^2 - ||x - \bar x||^2) + (||y - \bar x||^2 -
||y - \bar y||^2) \ge L^{-1} \delta_B (||\bar x - \bar y||/C_7),
\end{equation}
$$ C_7 = 2 {max} \lbrace 1, ||x - \bar y||, ||y - \bar x|| \rbrace .$$
It follows from (\ref{b12}) and (\ref{b13}) that
$$L^{-1} \delta_B (||\bar x - \bar y||/C_7) \leq  2C_4 \rho _B {(||x - y||)}.
$$
Finally, we have (\ref{b14}).  Theorem is proved.

Next we formulate the statement corresponding to the estimate of duality
mapping (\ref{d1}).
\begin {th} \label{th2}
In uniformly convex and uniformly smooth Banach space $B$ the following
estimate
\begin{equation} \label{b16}
||P_\Omega x - P_\Omega y|| \leq C \delta^{-1}_B (\rho _B {(8LC||x - y||)}),
\end{equation}
is satisfied where
$$C = 2 {max} \lbrace 1, ||x - \bar y||, ||y - \bar x|| \rbrace$$
\end{th}

We omit the proof of this Theorem because it is similar to the proof of
previous Theorem \ref{th1}.
\begin{re} For the Hilbert space $H$ one can write (\ref{b16}) in form
$$||P_\Omega x - P_\Omega y|| \leq 16LC^2||x - y||,$$
because  $\delta^{-1}_B (\cdot)$ and  $\rho _B (\cdot)$ are increasing
functions, $\rho _H (\tau) \leq \tau ^2 / 2$ and
$$ \epsilon ^2 /8  \leq \delta _{H}(\epsilon ) \leq \epsilon ^2 /4.$$
\end{re}
\begin{re}  It follows from (\ref {b14}) and (\ref {b16})  that  the
projection  operator $P_{\Omega }$ is
uniformly continuous on every bounded set of Banach space $B$.
\end{re}

It is interesting to compare the Bjornestal's estimate and a local version of
our estimate (\ref{b16}). In virtue of small distance
between $x$ and $y$ in (\ref{b2})  and condition $||x - \bar x|| = 1$ and
$||y - \bar y|| = 1$  one can consider that $C = 2$
in (\ref{b16}). Then
\begin{equation} \label{b15}
||P_\Omega x - P_\Omega y|| \leq 2 \delta^{-1}_B (\rho _B {(8LC||x - y||)}),
\;\; 1 < L < 3.18 .
\end{equation}
The constant in the brackets of
(\ref{b15}) is greater than the one in (\ref{b2}). It is natural
because (\ref{b16}) is global  estimate. Besides, the constants in (\ref{b2})
have been obtained  for the case $\Omega = M$, were $M$ is linear subspace of
$B$. (One might note that the constants in the estimates
(\ref{b2}), (\ref{b14}), (\ref{b16})) and (\ref{b15}) as a rule do not
have important meaning. There do have the orders of estimates,
and they are the same there).

\end{document}